\documentclass [11pt]{amsart} 
\usepackage{epic,eepic,latexsym, amssymb, amscd, amsfonts}

\input xy
\xyoption {all}

\def \uv {{\mathfrak M_{v}}}
\def \mv {{\mathsf M^{+}_{v}}}  
\def \uw {{\mathfrak M_{w}}}
\def \mw {{\mathsf M^{+}_{w}}}
\def \tv {{\Theta_{v}}}
\def \tw {{\Theta_{w}}}
\def \mvv {{\mathsf M^{-}_{v}}}
\def \mww {{\mathsf M^{-}_{w}}}
\def \adual {{\widehat A}}
\def \fm {{{\bf R}{\mathcal S}}}

\def \r {{{\bf R}}}

\newtheorem{theorem}{Theorem} 
\newtheorem {lemma}{Lemma}     
\newtheorem{conjecture}{Conjecture} 
\newtheorem {corollary}{Corollary} 
\newtheorem {proposition}{Proposition}
\newtheorem {fact}{Fact}
\newtheorem {assumption}{Assumption}

\theoremstyle{definition}
\newtheorem{remark}{Remark} 

\theoremstyle {definition}

\begin{document}

\title[Sheaves on abelian surfaces and Strange Duality]{Sheaves on abelian surfaces and 
Strange Duality}
\author {Alina Marian}
\address {School of Mathematics}
\address {Institute of Advanced Study}
\email {marian@math.ias.edu}
\author {Dragos Oprea}
\address {Department of Mathematics}
\address {Stanford University}
\email {oprea@math.stanford.edu}
\date{} 

\begin {abstract} We formulate three versions of a {\it strange duality} conjecture for sections of the Theta bundles on the moduli
spaces of sheaves on abelian surfaces. As supporting evidence, we check the equality of dimensions on dual moduli spaces, answering a
question raised by G\"{o}ttsche-Nakajima-Yoshioka \cite {GNY}.  
\end{abstract} \maketitle

\section{Introduction}

Let $(A, H)$ be a polarized abelian surface. In this paper, we consider the moduli spaces of Gieseker $H$-semistable sheaves on $A$, and
sections of the Theta line bundles defined over them.  

It will be convenient to bookkeep coherent sheaves $E$ on $A$ by their Mukai
vectors, setting $$v(E)=r +c_{1}(E)+\chi(E)\,\omega\in H^{2\star}(A),$$ where $\omega$ stands for the class of a point. As customary, we will equip the even cohomology $H^{2\star}(A)$ with the Mukai pairing. For any two vectors $x=(x_{0}, x_{2}, x_{4})\in
H^{2\star}(A)$ and $y=(y_{0}, y_{2}, y_{4})\in H^{2\star}(A)$, we set $$\langle x, y\rangle=-\int_{A} x^{\vee}\cup y=\int_{A}
\left(x_{2}y_{2}-x_{0}y_{4}-x_{4}y_{0}\right).$$ It follows that for any two sheaves $E$ and $F$, we have $$\chi(E, F)=\sum_{i=0}^{2} (-1)^{i}\text {
Ext}^{i}(E, F)=-\langle v(E), v(F)\rangle.$$

For an arbitrary $v\in
H^{2\star}(A)$, let us write $\uv$ for the moduli space of Gieseker $H$-semistable sheaves $E$ on $A$, with Mukai vector $v$. To keep
things simple, we will make the following assumption throughout: \begin{assumption}
\label {ass1}
\begin {itemize} \item [(i)] The polarization $H$ is generic i.e., it belongs to the complement of a locally finite union of hyperplane walls in the ample cone of $A$;\vskip.05in
\item [(ii)] The vector $v$ is a 
primitive  element of the lattice
$H^{2{\star}} (A, {\mathbb Z})$;\vskip.05in
\item [(iii)] The vector $v$ is
positive i.e., one of the following is true: \begin{itemize} \item $\text {rank }(v)>0$; \item $\text { rank 
}(v)=0
\text { and } c_{1}(v) \text { is effective, } \chi(v)\neq 0, \text { and } \langle v, v\rangle \neq 0, 4. $ 
\end {itemize} \end {itemize} \end {assumption}
The choice of generic polarization will play only a minor role in what follows, and as such, we will supress it from the notation. Note that $\text {(i)}$ and $\text {(ii)}$ together imply that $\uv$ is a smooth manifold consisting of stable sheaves only. Its dimension equals $2d_v + 2$, with \begin{equation}\label{dv}d_{v} =\frac{1}{2} \langle v, v\rangle.\end {equation}

The main characters of our story will be a collection of naturally defined Theta line bundles on $\uv$.  Consider the subgroup $v^{\perp}$
inside the holomorphic $K$-theory of $A$ generated by the sheaves $F$ whose Mukai vectors $w$ are orthogonal to $v$:
\begin{equation}\chi(v\otimes w)=-\langle v^{\vee}, w\rangle = 0.\label{ort}\end{equation} There is a morphism $$\Theta:{v}^{\perp}\to
\text {Pic}(\uv),\,\, \,\,\left[F\right]\to \Theta_F,$$ constructed and studied by Li and Le Potier \cite {Li} \cite {LP} in the case of surfaces, and 
also
by Dr\'{e}zet-Narasimhan in the case of curves \cite {dn}. The construction is easiest to explain assuming that $\uv$ is a fine moduli 
space,
such that $\mathcal E$ is the universal sheaf on $\uv \times A$. In this case, for a sheaf $F$ with Mukai vector $w$, we set \begin
{equation}\label{thw}\Theta_{F}=\det\r p_{!}(\mathcal E \otimes q^{\star} F)^{-1},\end {equation} where $p$ and $q$ are the two projections
from $\uv \times A$.  The orthogonality condition $\eqref{ort}$ is used to obtain a well defined line bundle $\Theta_{F}$, even in the
absence of universal structures, by descent from the Quot scheme.  Even though the line bundle $\Theta_{F}$ depends on the $K$-theory class
of $F$, the Chern class $c_{1}(\Theta_{F})$ depends only on the Mukai vector $w$. For simplicity of notation, in all cohomological
computations below, we will write $\Theta_{w}$ for any one of the line bundles $\Theta_{F}$ as above. The finer dependence of the line
bundles $\Theta_{F}$ on the sheaf $F$ will be discussed in more detail in Subsection \ref {thetas}.
\vskip.1in

Our goal in this paper is to compute the Euler characteristics of the line bundles $\Theta_{w}$, which can be interpreted as the $K$-theoretic Donaldson invariants of the abelian surface $A$. We will provide a simple expression for these Euler characteristics, valid in any rank.
We will thus answer a question raised in \cite {GNY}, as part of a general study of the rank two $K$-theoretic Donaldson invariants of
surfaces.

To explain the results, let us first fix a reference line bundle $\Lambda$ on $A$ with $c_{1}(\Lambda)=c_{1}(v)$. Then, we have a
well-defined determinant morphism \begin{equation}\label{det} \alpha^{+}_{\Lambda}: \uv\to \widehat A=\text {Pic}^{0}(A), \,\,\, E\to \det
E\otimes \Lambda^{-1}.\end{equation} Its fiber over the origin is the moduli space $\mathsf M_{v}^{+} (\Lambda)$ of sheaves with fixed
determinant $\Lambda$. The choice of the determinant is unimportant for our arguments, and therefore we will omit it from our notation when
no confusion is likely to arise. We will show:

\begin {theorem}\label{main} For any vectors $v$ and $w$ satisfying Assumption \ref{ass1}, and such that $\chi(v \otimes
w)=0$, we have \begin{equation} \label{abelianver} \chi (\mv, \Theta_{w})=\chi(\mw, \Theta_{v})=\frac{1}{2}\frac{c_{1}(v\otimes
w)^{2}}{d_{v}+d_{w}}\, \binom{d_{v}+d_{w}}{d_{v}}. \end{equation} \end {theorem}

There is yet another moduli space of interest to us, which is Fourier-Mukai `dual' to the one considered above. Letting $\mathcal P$ denote
the normalized Poincar\'e bundle on $A\times \widehat A$, the Fourier-Mukai transform is defined by \begin{equation}\label{fmtdef}\fm (E)=\r
p_{!}(\mathcal P\otimes q^{\star}E) \in {\bf D}(\widehat A),\end{equation} where $p, q$ are the two projections. With this understood, let
us set $$\alpha_{\Lambda}^{-}: \uv \to A, \,\,\, E\to \det \fm (E)\otimes \det \fm(\Lambda)^{-1}.$$ The fiber of the morphism 
$\alpha_{\Lambda}^{-}$
over the origin is denoted by $\mvv$, and parametrizes sheaves $E$ with fixed determinant of the Fourier-Mukai transform. We will prove:

\begin {theorem}\label{two}
For any vectors $v$ and $w$ as in Theorem \ref{main}, we have \begin{equation}\label{dual1}
\chi(\mvv, \Theta_{w})=\chi(\mww, \Theta_{v})=\frac{1}{2}\frac{c_{1}(\hat v\otimes \hat w)^{2}}{d_{v}+d_{w}}\binom{d_{v}+d_{w}}{d_{v}}. 
\end {equation}Here $\hat v$ and $\hat w$ denote the Fourier-Mukai transforms of the two vectors $v$ and $w$. 
\end {theorem}

Finally, we may consider the morphism $$\mathsf a_{v}=(\alpha^{+}_{\Lambda}, \alpha^{-}_{\Lambda}):\uv \to \widehat A\times A.$$ 
This is the Albanese map of the moduli space $\uv$, cf. \cite{Y}. Its fiber over the origin will henceforth be denoted by
$K_{v}$. We will show:

\begin{theorem}\label{three} Assume that the N\'{e}ron-Severi group of $A$ has rank $1$. With the same hypotheses as in Theorem \ref{main}, 
we
have \begin{equation}\label{kvmw}\chi(K_{v}, \tw)=\chi(\uw, \tv)=\frac{d_{v}^{2}}{d_{v}+d_{w}}\binom{d_{v}+d_{w}}{d_{v}}. \end {equation}
\end {theorem}

The manifest symmetry of the formulas in Theorems \ref{main}, \ref{two} and \ref {three} matches first of all that of their counterpart 
for 
the case of sheaves on a $K3$ surface. Indeed, the theta Euler characteristics for the moduli space of sheaves on a $K3$ were shown to be \cite {GNY}\cite {OG} 
\begin{equation}
\label{k3ver}
\chi(\mathfrak M_{v}, \Theta_{w})=\chi (\mathfrak M_{w}, \Theta_{v})=\binom{d_{v}+d_{w}+2}{d_{v}+1}.
\end{equation}

This symmetry further suggests a general {\it strange duality} for surfaces, reminiscent of the case of moduli spaces of bundles 
on curves. There, the analogous invariance of the Verlinde formula reflects a geometric isomorphism between generalized theta functions with 
dual ranks and levels \cite {MO} \cite{B}. In the case of sheaves on an abelian or $K3$ surface, it is tempting to assert, similarly, that whenever defined and nonzero, the morphisms
$$\mathsf{SD}^{\pm}: H^{0}(\mathsf M_{v}^{\pm},\tw)^{\vee}\to H^{0}(\mathsf M_{w}^{\pm}, \tv)$$ are isomorphisms. The same considerations
should apply to the companion morphism $$\mathsf {SD}:H^{0}(K_{v}, \tw)^{\vee}\to H^{0}(\uw, \tv).$$ In the above, for each of the three
pairs of moduli spaces {\it i.e.}, $(\mathsf M^{\pm}_{v}, \mathsf M_{w}^{\pm})$ and $(K_{v}, \uw)$, the line bundle $\tv$ stands for any one of the $\Theta_{E}$'s, for $E$ in the
corresponding moduli space of sheaves with Mukai vector $v$; similarly, $\tw$ is any one of the line bundles $\Theta_{F}$, for $F$ in the
dual moduli space of sheaves with Mukai vector $w$. We will review the definition of the three {\it strange duality} morphisms in Section
\ref{thetas} below, assuming that \begin {assumption}\label{ass4} For any two (semi)-stable sheaves $E$ and $F$ with Mukai vectors $v$ and
$w$, we have $$H^{2}(E\otimes F)=0.$$ This is automatic if $c_{1}(v\otimes w)\cdot H>0,$ by Serre duality and stability. \end {assumption}

In order to use the numerics provided by Theorems \ref {main}, \ref {two} and \ref {three}, one has to assume in addition that
the line bundle $\Theta_{w}$ has no higher cohomology on the various moduli spaces considered. The vanishing of higher cohomology is a
delicate question, which can be answered satisfactorily only in few cases. For smooth moduli spaces, or for moduli spaces with mild
singularities - {\it e.g.} rational - one may invoke standard vanishing theorems. These require the understanding of the positivity
properties of $\Theta_{w}$ {\it i.e.}, determining whether $\Theta_{w}$ is big and nef. In the case under study, smoothness is assumed,
bigness is easy to detect, and nefness is hoped for. This last point of nefness appears to be a subtle issue, even though the presence of the
holomorphic symplectic structure on the moduli spaces considered here makes the question more tractable. A study for the Hilbert scheme of two points on $K3$ surfaces and their deformations can be found for instance in \cite {HT}. Nevertheless, Le Potier \cite {LP} and Li \cite {Li} proved the following
results, which can be viewed as a higher dimensional generalization of the ampleness of the determinant line bundle on the moduli space of
bundles over a curve.

\begin {fact}
\begin{itemize}
\item [(i)] If $w$ has positive rank, and $c_{1}(w)$ is a high multiple of the polarization $H$, then $\Theta_{w}$ is relatively ample on the fibers of the determinant map $\alpha^{+}$. \vskip.07in

\item [(ii)] If $w$ has rank $0$, and $c_{1}(w)$ is a positive multiple of the polarization, then $\Theta_{w}$ is big and nef on the fibers of $\alpha^{+}$. \end{itemize}
\end {fact}

Similar results should hold for the morphism $\alpha^{-}$. When the Picard rank of $A$ is $1$, this is obtained for free in many cases, by the remarks following Conjecture \ref{rem1}\text {(ii)}.  

One may then speculate

\begin {conjecture} When Assumptions \ref{ass1} and \ref{ass4} are satisfied, the three morphisms $\mathsf {SD}^{+}, \mathsf {SD}^{-}$ and $\mathsf {SD}$ are either isomorphisms or zero. 
\end {conjecture}

This has the immediate 
\begin {corollary} As $E$ varies in $K_{v}$, the Theta sections $\Theta_{E}$ on the dual moduli space $\uw$ span the linear series $|\Theta_{v}|$. Same statements apply to the moduli spaces $\mw$ and $\mww$, letting $E$ vary in $\mv$ and $\mvv$ respectively. \end {corollary}

The Conjecture was demonstrated in a number of cases, in this and other geometric setups. An overview of the already existing arguments, as 
well as proofs of new cases, can be found in \cite{mo2}.

\subsection {Acknowledgements} 
The
calculations presented in this paper have as starting point Kota Yoshioka's extensive previous work on the subject. 

We would like to thank Jun Li for many conversations related to moduli spaces of sheaves on surfaces, and the NSF for financial support. A.M. is grateful to
Jun Li and the Stanford Mathematics Department for making possible a very pleasant stay at Stanford in the spring of 2007.

\section {Preliminaries on Theta divisors}

\subsection {The strange duality morphisms}\label {thetas} We review here the definition of the {\it strange duality} morphisms. We will fix
$\Lambda$ an arbitrary determinant, and we let $\widehat \Lambda = \det \fm (\Lambda).$ We recall the notation of the Introduction, letting
$\mv$ and $\mvv$ denote the moduli spaces of sheaves with determinant $\Lambda$ and determinant of the Fourier-Mukai transform equal to
$\widehat \Lambda$ respectively; $K_{v}$ consists of sheaves satisfying both requirements.

We explained in the introduction that the line bundle $\Theta_{F}$ only depends on the $K$-theory class of the reference sheaf $F$. 
We establish here the following more precise result.

\begin {lemma} \label{lem1} Consider a sheaf $F$ of Mukai vector $w$, and consider the line bundle $\Theta_{F}$ on the moduli space $\uv$. Then,

\begin {itemize}
\item [(i)] for $F\in \mw$, the restriction of $\Theta_{F}$ to $\mv$ is independent of the choice of $F$;
\item [(ii)] for $F\in \mww$, the restriction of $\Theta_{F}$ to $\mvv$ is independent of the choice of $F$;
\item [(iii)] for $F\in \uw$, the restriction of $\Theta_{F}$ to $K_{v}$ is independent of the choice of $F$.
\end {itemize}
\end {lemma}

\proof To prove $\text {(i)}$, pick two sheaves $F_{1}$ and $F_{2}$ with Mukai vector $w$ and the same determinant. Considering the virtual
element in $K$-theory $$\mathsf f=F_{1}-F_{2},$$ we need to show that $$\Theta_{F_{1}}\otimes \Theta_{F_{2}}^{-1}=\Theta_{\mathsf f}$$ is
trivial. We will show that in $K$-theory, \begin{equation}\label{fzw}\mathsf f=\mathcal O_{Z}-\mathcal O_{W}\end{equation} for two 
zero-dimensional schemes $Z$ and $W$, which necessarily have to be of the same length. This follows by induction on the rank of the $F$'s. The rank $0$ case is obvious. When the
rank is $1$, then $$F_{1}=L\otimes I_{Z}, F_{2}=L\otimes I_{W}$$ with $L=\det F_{1}=\det F_{2}$, and the result is immediate. For the inductive
step, note that it suffices to replace $F_{1}$ and $F_{2}$ by the twists $F_{1}(D)$ and $F_{2}(D)$, for some ample divisor $D$. In this case, we
reduce the rank by constructing exact sequences $$0\to \mathcal O_{A}\to F_{i}(D)\to F'_{i}\to 0$$ with $F'_{1},$ and $F'_{2}$ of the same lower
rank and the same determinant.  The claim then follows by the induction hypothesis applied to $F'_{1}-F'_{2}$. Once \eqref{fzw} is understood, it
suffices to assume that $Z$ and $W$ are supported on single points {\it i.e.,} that $\mathsf f$ is a formal sum $$\mathsf f=\sum_{z,w}\left
(\mathcal O_{z}-\mathcal O_{w}\right).$$

In this case, we will check $\Theta_{\mathsf f}$ is trivial by testing against any $S$-family $\mathcal E\to S\times A$ of sheaves with fixed
determinant $\Lambda$. In particular, the latter requirement implies that $$\det \mathcal E\cong M \boxtimes \Lambda$$ for some line bundle $M$ on $S$.
Then, the pullback of $\Theta_{\mathsf f}$ under the classifying morphism $S\to \mv$ is $$\det p_{!}(\mathcal E\otimes q^{\star}\mathsf
f)^{-1}=\left(\det \mathcal E\big{|}_{S\times \{z\}}\right)^{-1}\otimes\det \mathcal E\big{|}_{S\times \{w\}}\cong M^{-1}\otimes M\cong \mathcal
O_{S},$$ completing the proof of $\text {(i)}$.

For $\text{(ii)}$, the same reasoning applies to $$\fm ({\mathsf f}) =\fm (F_{1})-\fm(F_{2}),$$ which can be assumed to be the difference of
the structure sheaves of $z$ and $w$ on the dual abelian variety $\widehat A$. Therefore, $$\mathsf f=\mathcal P_{z}-\mathcal P_{w},$$ where
$\mathcal P_{z}$ and $\mathcal P_{w}$ are the line bundles on $A$ represented by $z$ and $w$. We need to check that $$\det p_{!}(\mathcal
E\otimes q^{\star}\mathcal P_{z})\cong \det p_{!}(\mathcal E\otimes q^{\star}P_{w}).$$ Consider the relative Fourier-Mukai sheaf on $S\times
\widehat A$ $$\det p_{13!}(p^{\star}_{12}\mathcal E\otimes p_{23}^{\star}\mathcal P)\cong M\boxtimes \widehat \Lambda,$$ for some line bundle
$M$ on $S$. Here the pushforward and pullbacks are taken along the projections from $S\times A\times \widehat A$. The conclusion follows by
looking at the isomorphic pullbacks of this sheaf over $S\times \{z\}$ and $S\times \{w\}$.

Finally, the third statement is obvious since $K_{v}$ is simply connected. Indeed, it suffices to check that $c_{1}(\Theta_{F})$ is independent of $F$. This
is a Grothendieck-Riemann-Roch computation, using the defining formula \eqref{thw}, with $\mathcal E$ replaced by a quasi universal sheaf, if needed.

\vskip.1in

\begin {remark}The Lemma above is sufficient to our purposes. It may be useful to have a more detailed understanding of how the Theta line bundles
$\Theta_{F}$ vary with $F$. For the case of curves, the requisite formulas were established by Dr\'{e}zet-Narasimhan \cite {dn}. We speculate that the following holds:
 
\begin {conjecture}\label{rem1} Consider two sheaves $F_{1}$ and $F_{2}$ with the same Mukai vector orthogonal to $v$. \begin{itemize} \item
[(i)] On $\mv$, we have $$\Theta_{F_{1}}=\Theta_{F_{2}}\otimes \left(\alpha^{-}\right)^{\star}\left (\det F_{1}\otimes \det F_{2}^{-1}\right).$$
\item [(ii)] On $\mvv$, we have $$\Theta_{F_{1}}=\Theta_{F_{2}}\otimes \left((-1)\circ \alpha^{+}\right)^{\star}\left(\det \fm (F_{1})\otimes
\det \fm (F_{2})^{-1}\right).$$ \item [(iii)] If $c_{1}(v)=0$, then on $\uv$ we have $$\Theta_{F_{1}}=\Theta_{F_{2}}\otimes \left((-1)\circ
\alpha^{+}\right)^{\star}\left(\det \fm (F_{1})\otimes \det \fm (F_{2})^{-1}\right)\otimes \left(\alpha^{-}\right)^{\star}\left (\det
F_{1}\otimes \det F_{2}^{-1}\right).$$ \end {itemize} \end {conjecture} Formula $\text {(i)}$ is easily checked on the Hilbert scheme of points
using that $$\Theta_{F}=(\alpha^{-})^{\star}(\det F)\otimes E^{\text { rank} F},$$ where $E$ is the exceptional divisor \cite {EGL}. Assuming
$\text {(i)}$, evidence for $\text{(ii)}$ is provided by the change of the Theta line bundles under Fourier-Mukai transform. Indeed, when the
Picard number of $A$ is $1$, Yoshioka \cite{Y} exhibited very general examples of birational isomorphisms between the moduli spaces $\mathsf M_{v}^{\pm}$ and
$\mathsf M_{\widehat v}^{\mp}$ on $A$ and $\widehat A$, interchanging the maps $\alpha^{+}$ and $\alpha^{-}$; arguments of Maciocia \cite{mac} can be used to show
that under this isomorphism the line bundle $\Theta_{F}$ corresponds to $\Theta_{(-1)^{\star}\widehat F}$, at least for generic $F$ satisfying
$WIT$.  Finally item $\text {(iii)}$ is consistent with $\text {(i)}$ and $\text {(ii)}$, and with Grothendieck-Riemann-Roch. It may be possible
to prove all three formulas using suitable degeneration arguments.  \vskip.1in \end {remark}

Assuming Lemma \ref {lem1}, it is now standard to define the three {\it strange duality} morphisms. The construction is contained in \cite {D} and \cite
{OG}, but we will review it briefly here for the sake of completeness.

Recall that for any pair of sheaves $(E, F)\in \uv \times \uw$ we have $$\chi(E\otimes F)=0.$$ We will assume furthermore that
\addtocounter{assumption}{-1} 
\begin {assumption} (modified version) \label{ass3}
\begin{itemize}
\item [(a)] either $H^{2}(E\otimes F)=0$; by stability this happens  if $c_{1}(E\otimes F).H>0$; 
\item [(b)] or $H^{0}(E \otimes F)=0$; by stability this happens if $c_{1}(E\otimes F).H<0$. 
\end {itemize}
\end {assumption}

To treat all cases at once, let us denote by $\mathcal M_{v}$ and $\mathcal M_{w}$ any one of the three pairs
$(\mv, \mw)$, $(\mvv, \mww)$ and $(K_{v}, \uw)$.  We set \begin{equation}\label{ll} \mathcal L_{w}=\begin{cases} \Theta_{F}, \text { for } F \in
\mathcal M_{w}, \text { if {\it Assumption \ref{ass3} (a)} holds},\\ \Theta_{F}^{-1}, \text { for } F \in \mathcal M_{w}, \text { if {\it
Assumption \ref {ass3} (b)} holds. }\end {cases}\end {equation} By Lemma \ref{lem1}, this is a well defined line bundle on $\mathcal M_{v}$. We
similarly define the line bundle $\mathcal L_{v}$ on $\mathcal M_{w}$. 

Descent arguments, presented in detail in D\v{a}nil\v{a}'s paper \cite {D}, show the existence of a natural divisor $$\Delta_{v, w}\hookrightarrow 
\mathcal
M_{v}\times \mathcal M_{w}$$ which is supported set-theoretically on the locus $$\Delta_{v,w}=\left \{(E, F)\in \mathcal M_{v}\times \mathcal M_{w}, \text
{such that } h^{1}(E\otimes F)\neq 0\right \}.$$ This divisor is obtained as the vanishing locus of a section of a naturally defined line bundle $\Theta_{v,
w}$ on $\mathcal M_{v}\times \mathcal M_{w}$. The splitting $$\Theta_{v,w}=\mathcal L_{w}\boxtimes\mathcal L_{v}$$ follows from Lemma \ref{lem1} by the
see-saw theorem. Therefore, $\Delta_{v,w}$ becomes an element of $$H^{0}(\mathcal M_{v}, \mathcal L_{w})\otimes H^{0}(\mathcal M_{w}, \mathcal L_{v})$$
inducing the duality morphism $$H^{0}(\mathcal M_{v}, \mathcal L_{w})^{\vee}\to H^{0}(\mathcal M_{w}, \mathcal L_{v}).$$ Note that when {\it Assumption \ref{ass3} (a)} holds, the construction gives rise to the three duality morphisms of the Introduction:\begin{equation}\label{sdpm}\mathsf {SD}^{\pm}:H^{0}(\mathsf M_{v}^{\pm}, \Theta_{w})^{\vee}\to H^{0}(\mathsf
M_{w}^{\pm}, \Theta_{v}), \text { and }\end{equation} \begin{equation}\label{sdn}\mathsf {SD}: H^{0}(K_{v}, \Theta_{w})^{\vee}\to H^{0}(\uw,
\Theta_{v}).\end {equation}
\section {Euler characteristics on Albanese fibers} 

\subsection {The Albanese map} In this section, we will compute the Euler characteristics of line bundles on $K_{v}$. We begin by reviewing a few facts about the Albanese map of $\uv$. Recall from the introduction the modified 
determinant
morphism $$\alpha^{+}_{\Lambda}:\uv \to \adual, E\to \det E \otimes \Lambda^{-1},$$ and its Fourier-Mukai 'dual' $$\alpha^{-}_{\Lambda}:
\uv \to A, \,\,\, E\to \det \fm (E)\otimes \widehat \Lambda^{-1}.$$ Putting these two morphisms together, we obtain the map
\begin{equation}\label{alb}\mathsf a_{v}=(\alpha^{+}, \alpha^{-}):\uv \to \widehat A\times A.\end {equation} Yoshioka
proved that $\mathsf a_{v}$ is the Albanese map of the moduli space $\uv$ \cite {Y}.

The morphism \eqref{alb} is easiest to understand for the vector $v=(1,0, n)$. Then, the moduli space $\uv$ is isomorphic to the product $\widehat A\times A^{[n]}$ of the dual abelian variety and the Hilbert scheme $A^{[n]}$ of points on $A$. The 
morphism $\mathsf a_{v}$ can be identified with $$1\times s: \widehat A \times A^{[n]} \to \widehat A\times A$$ where the first 
map is the identity, while the second is induced by summation on the abelian surface. That is, for a zero-cycle $Z$ supported on points $z_{i}$ with length $n_{i}$, we let $$s:A^{[n]}\to A,\,\,\, s([Z])=\sum_{i} n_{i}z_{i}.$$  The fiber of $s$ over $(0,0)$ is the generalized Kummer variety $K_{n-1}$ of dimension $n-1$. 

In general, Yoshioka studied the fiber of the Albanese map $\mathsf a_{v}$ over the origin $$K_{v}=\mathsf a_{v}^{-1}(0, 0),$$ under the
assumption that the vector $v$ is primitive and positive, in the sense discussed in the Introduction. In this situation, and when
$d_{v}\geq 3$, Yoshioka proved that \begin {itemize} \item $K_{v}$ is an irreducible holomorphic symplectic manifold, deformation equivalent to the
generalized Kummer surface $K_{n-1}$, with $n=\langle v, v \rangle/2$. \item There is an isomorphism $$H^{2}(K_{v}, \mathbb Z)\cong \text
{Pic} (K_{v}) \cong v^{\perp}.$$ The vector $w$ in $v^{\perp}$ corresponds to the line bundle $\Theta_{w}$ restricted to $K_{v}$.  \item
Moreover, if one endows $H^{2}(K_{v},\mathbb Z)$ with the Beauville-Bogomolov form, and $v^{\perp}$ with the intersection form, the above isomorphism
is an isometry.  \end {itemize}

\subsection {Generalized Kummer varieties.} We start the calculation of Euler characteristics by considering the case of line bundles on the generalized
Kummer varieties. In other words, we assume that $v=(1,0,n)$. 

For any divisor $D$ on the abelian surface $A$, we let $D_{(n)}$ be the divisor on $A^{[n]}$ consisting
of zero-cycles which intersect $D$. This divisor is a pull-back under the support morphism 
$$f: A^{[n]}\to A^{(n)},$$
from the symmetric power $A^{(n)}$ of $A$, 
$$D_{(n)}=f^{\star} \left(D\boxtimes D \boxtimes \ldots \boxtimes D\right)^{S_{n}}.$$ Further, let $E$ be
the exceptional divisor of $A^{[n]}$ consisting of schemes with two coincident points in their support. Any line bundle on the Hilbert
scheme is of the form $D_{(n)}\otimes E^{r}.$ We will denote by the same symbol the restriction of these bundles to the generalized Kummer variety. We also do not distinguish notationally between line bundles and divisors. 

\begin {lemma}
\label{first0}
$$\chi(K_{n-1}, D_{(n)}\otimes E^{r})=n\binom{\chi(D)-(r^{2}-1)n-1}{n-1}$$
\end {lemma}

\proof The expression given by the Lemma is a consequence of the known formula 
\begin{equation}\label{pull2}\chi (A^{[n]}, D_{(n)}\otimes E^{r})=\frac{\chi(D)}{n} \binom {\chi(D)-(r^{2}-1)n -1}{n-1},\end{equation}
which was deduced in \cite{EGL}. To relate the two, we use the cartesian diagram

\begin{center}
$\xymatrix{K_{n-1}\times A \ar[r]^{\sigma} \ar[d]^{p} & A^{[n]} \ar[d]^{s}
\\ A \ar[r]^{n} &  A}.$
\end{center}
The upper horizontal map is $$\sigma:K_{n-1}\times A\to A^{[n]}, \,\, (Z, a) \mapsto t_{a}^{\star} Z,$$ while the bottom morphism is the 
multiplication by $n$ in the abelian surface. It follows that $\sigma$ is an \'{e}tale cover of degree $n^{4}$. 

By the see-saw theorem, we find $$\sigma^{\star} D_{(n)}=D_{(n)}\boxtimes D^{\otimes n}, \text { while } \sigma^{\star}E=E\boxtimes 
\mathcal O_{A}.$$
Therefore \begin{equation}\label{pull1} \chi(A^{[n]}, D_{(n)}\otimes E^{r})=\frac{1}{n^{4}} \chi(K_{n-1}, D_{(n)}\otimes E^{r}) \,\chi(A,
D^{n})=\frac{\chi(D)}{n^{2}} \chi(K_{n-1}, D_{(n)}\otimes E^{r}).\end {equation} 
Putting \eqref
{pull2} and \eqref{pull1} together we obtain the Lemma.\vskip.1in

\subsection {General Albanese fibers.} We can now consider the case of an arbitrary vector $v$. We claim 

\begin{proposition}\label{kv} If $d_{v}\neq 0$, then 
$$\chi (K_{v}, \Theta_{w})=\frac{d_{v}^{2}}{d_{v}+d_{w}} \binom{d_{v}+d_{w}}{d_{v}}.$$
\end {proposition}

\proof We prove here the statement when $d_{v}\neq 2$. The case $d_v = 2$ will be considered separately in the next subsection. 

When $d_{v}=1$ the Proposition follows immediately from Mukai and Yoshioka's results \cite{Muk}, \cite{Y3}, as they have proved that the Albanese 
map $\mathsf 
a_{v}:\uv\to \widehat A\times A$ is an isomorphism. Then $K_{v}$ is a point, and both sides of the equation in Proposition \ref{kv} equal $1$. 

When $d_{v}\geq 3$, we follow the same arguments as in the case of $K3$ surfaces. We will make use of the Beauville-Bogomolov form $B$. This
quadratic form is defined on the second cohomology of any irreducible holomorphic symplectic manifold, and can be considered as a generalization of the
intersection pairing on $K3$ surfaces. In the case of the generalized Kummer varieties $K_{n-1}$, the form $ B$ gives an orthogonal
decomposition $$H^{2}(K_{n-1}, \mathbb Z)= H^{2}(A, \mathbb Z) \oplus \mathbb Z[E]$$ such that $$B(c_1 (D_{(n)}))=D^{2}, \,\,\,\, 
B(c_1(E))=-2n.$$ In
particular, $$B(c_1(D_{(n)}\otimes E^{r}))=2\left(\chi(D)-r^{2}n\right).$$ Therefore, the result of Lemma \ref{first0} can be restated as
$$\chi(K_{n-1}, L)=n \binom {\frac{B(c_{1}(L))}{2}+n-1}{n-1},$$ for any line bundle $L=D_{(n)}\otimes E^{r}$ on the Kummer variety
$K_{n-1}$.

To get the result of the Proposition, we will use the fact that the Euler characteristics $\chi(X, \mathcal L)$ of any line bundle on an irreducible
holomorphic symplectic manifold $X$ can be expressed as a universal polynomial in the Beauville-Bogomolov form $B(c_1 (\mathcal L))$ \cite {H1} {\it i.e.,}
a polynomial depending only on the underlying holomorphic symplectic manifold. This polynomial is an invariant of the deformation type.
Moreover, the Beauville-Bogomolov form is also invariant under deformations. Now, by \cite{Y}, we know that $K_{v}$ is deformation equivalent to the
generalized Kummer variety $K_{n-1}$, for $n=\langle v, v\rangle/2$, and that $$B(c_{1}(\Theta_{w}))=\langle w, w\rangle.$$ Therefore,
$$\chi(K_{v}, \Theta_{w})=d_{v} \binom {\frac{B(c_{1}(\Theta_{w}))}{2} + d_{v}-1}{d_{v}-1}=d_{v}\binom
{d_{v}+d_{w}-1}{d_{v}-1}=\frac{d_{v}^{2}}{d_{v}+d_{w}}\binom {d_{v}+d_{w}}{d_{v}},$$ completing the proof of the Proposition when $d_v \neq 2.$

\subsection {Two-dimensional Albanese fibers} In this subsection, we establish Proposition \ref{kv} when  $d_{v} = 2,$ by analyzing the fairly involved special 
geometry of the 
situation. The case $r= 1$ is covered by Lemma \ref{first0}, so we assume that $r\geq 2$. 
$K_{v}$ is now a $K3$ surface. It suffices to show that $$\chi(K_{v}, \Theta_{w})=\frac{c_{1}(\Theta_{w})^{2}}{2}+2=2d_{w}+2$$ or
equivalently, \begin{equation}\label{dv2}c_{1}(\Theta_{w})^{2}=2\langle w, w\rangle.\end{equation}

Yoshioka identifies the $K3$ surface $K_{v}$ as a Fourier-Mukai partner of the Kummer surface $X$ associated to $A$ \cite {Y}. We 
review his construction below. To fix the notation, consider the following diagram \begin {center} $\xymatrix{ \bigsqcup_{i=1}^{16}
E_{i}\ar@{^{(}->}[r]& \tilde A\ar[r]^{p}\ar[d]^{\pi} & X& \bigsqcup_{i=1}^{16}C_{i}\ar@{_{(}->}[l]\\ & A}$ \end {center} Here $\tilde A$ is
the blowup of $A$ at the $16$ two-torsion points, with exceptional divisors denoted by $E_{i}$.  Let $$j:\bigsqcup_{i=1}^{16}E_{i}\to
\tilde A$$ denote the inclusion of all $16$ exceptional divisors in $\tilde A$. Finally, $p$ is the morphism quotienting the $\mathbb 
Z/2\mathbb
Z$ automorphisms, and $C_{i}$ are the $(-2)$ curves on $X$ which are the images of the exceptional divisors $E_{i}$ under $p$.

Yoshioka proves that for a suitable isotropic Mukai vector $\tau$ on $X$ {\it i.e.}, $\langle \tau, \tau
\rangle = 0$, and a suitable polarization $L$, the following isomorphism holds \begin{equation}\label{isoyo}K_{v}\cong \mathfrak M_{X}(\tau).\end {equation} Here $\mathfrak M_{X}(\tau)$ denotes the moduli space of $L$-semistable sheaves
on $X$, with Mukai vector $\tau$. We will use Yoshioka's explicit isomorphism to identify the theta bundle $\Theta_{w}$ on the moduli space $\mathfrak M_{X}(\tau)$. 

We will explain the argument when $r+ c_{1}(v)$ is indivisible. The remaining case when $r+c_{1}(v)$ equals twice a primitive class is entirely similar. Let us assume first that $r\neq 2$, and that either $r$ and $\chi(v)$ are both even, or $r$ is odd. The isomorphism \eqref{isoyo} associates to each $F \in \mathfrak M_{X}(\tau)$ a sheaf $E$ on $A$, via elementary modifications along the exceptional divisors on the blowup $\widetilde {A}$. Concretely, the sheaf $p^{\star}F|_{E_{i}}$ splits as a sum \begin{equation}\label{yo2}p^{\star}F|_{E_{i}}\cong \mathcal 
O_{E_{i}}(-1)^{\oplus{a_{i}}}\oplus
\mathcal O_{E_{i}}^{\oplus(r-a_{i})},\end{equation} for suitable integers $a_{i}.$ 
Then, $E$ is defined by the exact sequence
\begin{equation}\label{yo}0\to \pi^{\star}E\to p^{\star}F\to j_{\star}\left(\bigoplus_{i=1}^{16} \mathcal 
O_{E_{i}}(-1)^{\oplus a_{i}}\right)\to
0.\end{equation} The assignment $$\mathfrak M_{X}(\tau)\ni F\to E\in \uv,$$ establishes an isomorphism onto the image $K_{v}$. 

In fact, Yoshioka's construction works in families, giving a natural transformation between the moduli functors $$\underline{\mathfrak
M}_{X}(\tau)\to \underline{K}_{v}.$$ The exact description of this transformation will be useful later. In what follows, let us agree that
the base change of various morphisms to an arbitrary base $S$ will be decorated by overlines. Fix any flat $S$-family $\mathcal F$ of
sheaves in $\underline {\mathfrak M}_{X}(\tau) (S)$. Define the sheaf $\mathcal G$ on the union $\bigsqcup_{i=1}^{16} E_i \times S$ via the exact sequence 
\begin{equation}\label{yo3}0\to
\bar\pi^{\star}\bar\pi_{\star} \bar j^{\star}\bar p^{\star}\mathcal F\to \bar j^{\star}\bar p^{\star}\mathcal F\to \mathcal G\to 0.\end
{equation} The short exact sequence \begin{equation}\label{yo4}0\to \bar \pi^{\star}\mathcal E\to \bar p^{\star}\mathcal F\to \bar
j_{\star} \mathcal G\to 0\end{equation} then defines a new $S$-family $\mathcal E$ in $\underline {K}_{v}(S)$.

Now let \begin{equation}\label{zetadef}\zeta=\text {ch}\left(p_{!}\left(\pi^{\star} W\right)\right)(1+\omega)\in H^{\star}(X),\end{equation} be the Mukai vector of the pushforward $p_{!}(\pi^{\star} W)$, for an arbitrary sheaf $W$ on $A$ with Mukai vector $w$. Using the exact sequence \eqref{yo}, and the fact that
$$\chi(j_{\star}\mathcal O_{E_{i}}(-1)\otimes \pi^{\star} w)=0,$$ we find $$2 \, \chi(\tau\otimes \zeta)=\chi(p^{\star}\tau \otimes
\pi^{\star}w)=\chi(\pi^{\star}(v\otimes w))=\chi(v\otimes w)=0.$$

The above computation shows that $\zeta$ defines a Theta line bundle $\Theta_{\zeta}$ on $\mathfrak M_{X}(\tau).$ We claim that under the
isomorphism $K_{v} \cong \mathfrak M_{X}(\tau)$ we have an identification \begin{equation}\label{ref}\Theta_{w}\cong\Theta_{\zeta}.\end{equation} As for any good quotient, the Picard
group of the moduli scheme $\mathfrak M_{X}(\tau)$ injects into that of the moduli functor $\underline {\mathfrak M}_{X}(\tau)$. Therefore, it
suffices to check equality of the two line bundles $\Theta_{w}$ and $\Theta_{\zeta}$ over arbitrary base schemes $S$, and for arbitrary $S$-families
$\mathcal F$ of $ \underline{\mathfrak M}_{X}(\tau)$. 

Let $q:S\times \tilde A\to S$ and $\widetilde q:S\times X\to S$ be the two projections. The exact sequence \eqref{yo4} and the push-pull formula then
give $$\Theta_{w}=\det \r q_{!}\left(\bar\pi^{\star}\mathcal E\otimes \text {pr}_{A}^{\star} w\right)^{-1}=\det \r q_{!}\left(\bar
p^{\star}\mathcal F\otimes \text {pr}_{A}^{\star}w\right)^{-1}=\det \r \widetilde q_{!}\left(\mathcal F\otimes
\text{pr}^{\star}_{X}\zeta\right)^{-1}=\Theta_{\zeta}.$$ Here, we used that the contribution of the last term of \eqref{yo4} vanishes. Indeed, since $$\bar
j^{\star}\text{pr}_{A}^{\star}w=\text {rank }(w)\cdot 1,$$ we have $$\det \r q_{!} \left(\bar j_{\star} \mathcal G \otimes
\text{pr}_{A}^{\star}w\right)=\text {rank }(w)\cdot \det \r \left(q\bar j\right)_{!}(\mathcal G)=0.$$ The last equality follows from the fact
that all direct images of $\mathcal G$ vanish. This is implied by the base change theorem, observing that the restriction of $\mathcal G$ to
each fiber of the morphism $q\bar j:S\times E_{i}\to S$ splits as a sum of line bundles $\mathcal O_{E_{i}}(-1).$ In turn this latter fact is a
consequence of the defining exact sequence \eqref{yo4}, in conjunction with equation \eqref{yo2}. 

To complete the proof, recall that Mukai \cite{mukai3} established an isometric isomorphism $$H^{2}(\mathfrak M_{X}(\tau))\cong
\tau^{\perp}/\tau$$ where the left hand side is endowed with the intersection pairing, while the right hand side carries the Mukai form induced
from the cohomology $H^{\star}(X)$. Then, $$c_{1}(\Theta_{w})^{2}=c_{1}\left(\Theta_{\zeta}\right)^{2}=\langle \zeta, \zeta\rangle = 2 \langle
w, w \rangle.$$ This proves \eqref{dv2}.

The case $r\neq 2$ and $\chi(v)$ odd is entirely similar. In this case, the exact sequence \eqref{yo} is replaced by $$0\to \pi^{\star}E\to p^{\star}F \,(E_{1})\to j_{\star}\left(\bigoplus_{i=1}^{16} \mathcal 
O_{E_{i}}(-1)^{\oplus a_{i}}\right)\to
0.$$ The argument identifying the Theta bundles carries through, for the vector \begin{equation}\label{zetadef2}\zeta=\text {ch}(p_{!}(\pi^{\star} W(E_{1}))(1+\omega).\end{equation}

The case $r=2$ requires a different discussion, since in this case, the description of the isomorphism \eqref{isoyo} via the assignment $F\to E$ is valid only on the complement of four rational curves $R_{i},$ $1\leq i\leq 4$. In fact, one cannot pick an isotropic vector $\tau$ such that for each of the $16$ exceptional divisors, the 
splitting type \eqref{yo2} is independent of the choice of a point $[F] \in {\mathfrak M}_{X} (\tau).$    
At best, for a suitable $\tau$, the rigid splitting $$p^{\star}F|_{E_{i}}\cong \mathcal
O_{E_{i}}\oplus\mathcal O_{E_{i}}(-1)$$ holds for $12$ exceptional divisors $E_{i}, \, 5\leq i\leq 16.$ For the remaining four divisors
$E_{i}, 1\leq i\leq 4$, the splitting type varies within the moduli space. 

When $\chi(v)$ is even, for generic $F$ in $\mathfrak M_{X}(\tau)$, the splitting type
is \begin{equation} \label{nong2}p^{\star}F|_{E_{i}}=\mathcal O_{E_{i}}\oplus \mathcal O_{E_{i}},\,\, 1\leq i\leq 4.\end{equation} The loci
of non-generic splitting give rational curves $R_{i}$ in $K_{v}$. 
Indeed, the nongeneric splitting is \begin{equation}\label{nong}p^{\star}F|_{E_{i}}=\mathcal O_{E_{i}}(-1)\oplus \mathcal O_{E_{i}}(1).\end{equation}
These $F$'s are shown to sit in exact sequences \begin{equation}\label{yo7}0\to G_{i}\to F\to \mathcal O_{C_{i}}(-1)\to 0\end{equation} for
certain rigid stable bundles $G_{i}$ on the $K3$ surface $X$, cf. Lemma $4.23$ in \cite {Y}. The nontrivial extensions \eqref{yo7} are parametrized by a rational curve
$$R_{i}\cong \mathbb P(\text {Ext}^{1}_{X}(O_{C_{i}}(-1), G_{i})).$$ 

When $\chi(v)$ is odd, all these statements are true for the exceptional divisors $E_{2}, E_{3}, E_{4}$, but equations \eqref{nong2} and \eqref{nong} fail for $E_{1}$. In fact, generically \begin{equation}\label{nong3}p^{\star}F|_{E_{1}}=\mathcal O_{E_{1}}(1)\oplus \mathcal O_{E_{1}}(1),\end{equation} while nongenerically \begin{equation}\label{nong1}p^{\star}F|_{E_{1}}=\mathcal O_{E_{1}}\oplus\mathcal O_{E_{1}}(2).\end{equation} The nongeneric splitting occurs along the rational curve $$R_{1}=\mathbb P(\text {Ext}^{1}(\mathcal O_{C_{1}}, G_{1}))$$ parametrizing extensions of the type \begin{equation}\label{yo10}0\to G_{1}\to F\to \mathcal O_{C_{1}}\to 0,\end{equation} for some rigid vector bundle $G_{1}$ on $X$.

We claim that the Theta bundles agree in this case as well {\it i.e.}, we check that the isomorphism \eqref{ref} $$\Theta_{w}\cong \Theta_{\zeta}$$ is satisfied. Let us first discuss the case when $\chi(v)$ is even, with $\zeta$ given by \eqref{zetadef}. To begin, $\Theta_{w}$ and $\Theta_{\zeta}$ agree on the complement of the four rational curves $R_{i}$, $1\leq i\leq 4$, since the exact sequence \eqref{yo} is valid outside these curves. We will check that the Theta bundles agree along the curves $R_{i}$ as well. Precisely, we claim that \begin {equation}\label {first} c_{1}(\Theta_{w})\cdot R_{i}= c_{1}(\Theta_{\zeta})\cdot R_{i}=s,\,\,\, 1\leq i\leq 4,\end {equation} with $$s=\text {rank } w.$$ Moreover, the four curves $R_{i}, 1\leq i \leq 4$, are disjoint. These facts will establish the isomorphism \eqref{ref}.

To calculate the first intersection in \eqref{first}, we will use the explicit description of the rational curves $R_{i}$ in the moduli space $K_{v}$. Specifically, Yoshioka notes that the curve $R_{i}$ corresponds to those sheaves $E$ on $A$ which fail to be locally free at a two-torsion point $x_{i}$. We can construct these sheaves as elementary modifications of a fixed $V$: \begin{equation}\label{construct}0\to E\to V\to \mathcal O_{\{x_{i}\}}\to 0.\end{equation} Since the middle sheaf $V$ has Mukai vector $v+\omega$, it sits in a moduli space of dimension $2$. Mukai showed that all such $V$'s are locally free \cite {mukai3}. Therefore, $E$ is not locally free at $x_{i}$, but it is locally free elsewhere. These nonlocally free elementary modifications are moreover parametrized by a $\mathbb P^{1}$, which should therefore be the rational curve $R_{i}$ above.  Moreover, the argument shows that the four rational curves $R_{i}$, $1\leq i\leq 4$ are disjoint.

The universal structure on $R_{i}\times A$, associated to the elementary modifications \eqref{construct}, becomes $$0\to \mathcal E\to \text{pr}^{\star}_{A} V\to \mathcal O_{R_{i}}(1)\boxtimes \mathcal O_{\{x_{i}\}}\to 0.$$ Therefore, $$c_{1}(\Theta_{w})\cdot R_{i}=-c_{1}(p_{!}(\mathcal E\otimes \text{pr}^{\star}_{A} w))=c_{1}\left(p_{!}\left(\mathcal O_{R_{i}}(1)\boxtimes \left(\mathcal O_{\{x_{i}\}}\otimes  w\right)\right)\right)=c_{1}(\mathcal O_{R_{i}}(1)^{\oplus s})=s.$$ 
   
To prove the second equality of \eqref{first}, we will use the description of the rational curves $R_{i}$ provided by equation \eqref{yo7}. The universal extension $$0\to \text{pr}^{\star}_{X}G_{i}\to \mathcal F\to \mathcal O_{R_{i}}(-1)\boxtimes \mathcal O_{C_{i}}(-1)\to 0$$ on $R_{i}\times X$ restricts to \eqref{yo7} on the fibers of the projection $p: R_{i}\times X\to R_{i}$. Using this exact sequence, we compute $$c_{1}(\Theta_{\zeta})\cdot R_{i}=-c_{1}(p_{!}(\mathcal F\otimes \text {pr}^{\star}_{X} \zeta))=-c_{1}\left(p_{!}\left(\mathcal O_{R_{i}}(-1)\boxtimes \left(\mathcal O_{C_{i}}(-1)\otimes \zeta\right)\right)\right)$$ $$=-c_{1}(\mathcal O_{R_{i}}(-1))\,\chi(\mathcal O_{C_{i}}(-1)\otimes \zeta)=c_1(\zeta)\cdot C_{i}=s.$$ The last evaluation follows from \eqref{zetadef} via Riemann-Roch. 


When $\chi(v)$ is odd, the numerics are slightly different, but \eqref{ref} still holds for the vector $\zeta$ given by \eqref{zetadef2}. In this case, we check that $$c_{1}(\Theta_{\zeta})\cdot R_{1}=s,$$ using equation \eqref{yo10}, instead of \eqref{yo7}. 

This completes our analysis of the two-dimensional Albanese fibers, establishing Proposition \ref{kv}.

\section {Sheaves with fixed determinant} In this section we will prove Theorem \ref{main}. We begin by fixing the notation. Specifically, let us write $r,
\Lambda, \chi$ for the rank, determinant and Euler characteristic of the vector $v$. The notation $r', \Lambda', \chi'$ will be used for the vector $w$. The
orthogonality of $v$ and $w$ translates into \begin{equation}\label{relation}r'\chi+c_{1}(\Lambda)\cdot c_{1}(\Lambda')+r\chi'=0.\end {equation}
Let $\mathcal P$ be the normalized Poincar\'{e} bundle on $A\times \widehat A$. We make the convention that $x$ will stand for a point of $A$, while $y$ 
will
be a point of $\widehat A$. We will write $$\mathcal P_{x}=\mathcal P|_{\{x\}\times \widehat A},\,\,\,\mathcal P_{y}=\mathcal P|_{A\times \{y\}}\cong y.$$ We
denote by $t_{x}$ and $t_{y}$ the translations by $x$ and $y$ on the abelian varieties $A$ and $\widehat A$ respectively.

The following two facts about the Fourier-Mukai transform of an arbitrary $ E\in {\bf D} (A)$, proved in \cite {Muk}, will be used below:
\begin{equation}\label{fmpp1}\fm (E\otimes \mathcal P_{y})=t_{y}^{\star}\, \fm(E),\end{equation}
\begin{equation}\label{fmpp}\fm(t_{x}^{\star}E)=\fm(E)\otimes \mathcal P_{-x}.\end {equation}

It is moreover useful to recall that the two line bundles $\Lambda$ and $\widehat \Lambda$ standardly induce morphisms $$\Phi_{\Lambda}:A\to
\widehat A, \,\,x\mapsto t_{x}^{\star}\Lambda \otimes \Lambda^{-1}, \text { and }$$ $$\Phi_{\widehat \Lambda}:\widehat A\to A,\,\, y\to t_{y}^{\star}\widehat
\Lambda\otimes \widehat \Lambda^{-1},$$ satisfying \cite{Y} \begin{equation}\label{phis}\Phi_{\Lambda}\circ \Phi_{\widehat
\Lambda}=-\chi(\Lambda)1,\,\, \Phi_{\widehat \Lambda}\circ \Phi_{\Lambda}=-\chi(\Lambda) 1.\end{equation}

To start the proof of Theorem \ref{main}, consider the diagram  \begin{center} $\xymatrix{K_{v}\times A \ar[r]^{\Phi^{+}} \ar[d]^{p} & \mv
\ar[d]^{\alpha^{-}} \\ A \ar[r]^{\Psi^{+}} & A}.$\label{diag} \end{center} The upper horizontal map is given by $$\Phi^{+}(E,
x)=t_{rx}^{\star}E \otimes t_{x}^{\star}\Lambda^{-1}\otimes \Lambda.$$ This is well defined since $$\det \Phi^{+}(E, x)=t_{rx}^{\star}\Lambda
\otimes \left(t_{x}\Lambda^{-1}\otimes \Lambda\right)^{r}=\Lambda.$$

\begin{lemma}\label{multlemma} The morphism $\Psi^{+}$ is the multiplication by $d_{v}$ in the abelian variety. 
\end {lemma}

\proof Using \eqref{fmpp1} and \eqref {fmpp}, we compute \begin{eqnarray*}\alpha^{-} \circ \Phi^{+} (E, x)&=&\det \fm(t_{rx}^{\star}E\otimes
t_{x}^{\star}\Lambda^{-1}\otimes \Lambda)\otimes \widehat \Lambda^{-1}=\det \fm (t_{rx}^{\star}E \otimes \mathcal P_{-\Phi_{\Lambda}(x)})\otimes
\widehat\Lambda^{-1}\\ &=&\det \left(t_{-\Phi_{\Lambda}(x)}^{\star} \fm (t_{rx}^{\star}E)\right) \otimes \widehat \Lambda^{-1}=t_{\Phi_{\Lambda}(-x)} ^{\star}\det
\fm(t_{rx}^{\star}E) \otimes \widehat \Lambda^{-1}\\ &=&t_{\Phi_{\Lambda}(-x)} ^{\star}\det\left(\fm(E)\otimes \mathcal P_{-rx}\right) \otimes \widehat\Lambda^{-1}
\\&=&t_{\Phi_{\Lambda}(-x)}^{\star} \left(\det \fm(E)\otimes \mathcal P_{-rx}^{\chi}\right)\otimes \widehat \Lambda^{-1}\\ &=&t_{\Phi_{\Lambda}(-x)}^{\star}
\widehat \Lambda \otimes \mathcal P_{-r\chi x}\otimes \widehat \Lambda^{-1}=\Phi_{\widehat \Lambda}(\Phi_{\Lambda}(-x))\otimes \mathcal P_{-r\chi
x}\\&=&\chi(\Lambda) x\otimes \mathcal P_{-r\chi x}=(\chi(\Lambda)-r\chi) x=d_{v}x.\end{eqnarray*} The first equality on the penultimate line follows from the fact that the Poincar\'{e} bundle $\mathcal P_{x}$ is invariant under translations \cite {M}. Equation \eqref{phis} was used in the last line.

\begin {lemma} When $d_{v}\neq 0$, the diagram above is cartesian. Therefore, the morphism $\Phi^{+}$ has degree $d_{v}^{4}$.
\label{degree1}
\end {lemma} 

\proof This is almost immediate. Together, $\Phi^{+}$ and $p$ give rise to a morphism $i:K_{v}\times A\to \mv \times_{A, (\alpha^{-}, \Psi^{+})}
A.$ We show that $i$ is an isomorphism. Since $\Psi^{+}$ is \'{e}tale, the natural morphism $\mv \times_{A, (\alpha^{-}, \Psi^{+})} A\to \mv$ is also 
\'{e}tale, 
so
the fibered product $\mv \times_{A, (\alpha^{-}, \Psi^{+})} A$ is smooth. The fibered product is also connected, as it follows by looking at the connected fibers of the projection to $A$; note that the projection is surjective, as $\alpha^{-}$ has this property, according to the previous Lemma. Therefore, it suffices to check that $i$ is injective. If $i(E, x)=i(E', x')$ then,
by composing $i$ with $\Phi^{+}$ and $p$, we see that $$t_{rx}^{\star}E\otimes t_{x}^{\star} \Lambda^{-1}\otimes \Lambda=t_{rx'}^{\star}E'\otimes
t_{x'}^{\star}\Lambda^{-1} \otimes \Lambda, \text { and } x=x'.$$ This immediately implies $E=E'$ as well. The diagram is therefore cartesian. \vskip.1in

\begin {proposition} \label {split} We have $$\left(\Phi^{+}\right)^{\star}\Theta_{w}\cong \Theta_{w}\boxtimes \mathcal L^{+}$$ where $\mathcal L^{+}$ is a line 
bundle
on $A$ with $$c_{1}(\mathcal L^{+})=- d_{v}c_{1}(v\otimes w).$$ \end {proposition} \proof This follows by the see-saw theorem. Letting
$\Phi_{x}=\Phi^{+}|_{K_{v}\times \{x\}}$, we claim that the pullback $\Phi_{x}^{\star}\tw$ is independent of $x$, and therefore, by specializing to $x=0$, it
should coincide with $\tw$. Since $K_{v}$ is simply connected, it suffices to check that the Chern class $c_{1}(\Phi_{x}^{\star}\Theta_{w})$ is independent
of $x$. This is clear when a universal sheaf $\mathcal E$ exists on $\mv \times A.$ Indeed, for $F$ a sheaf on 
$A$ with Mukai vector $w$, $$\Phi_{x}^{\star}\Theta_{w}=\Phi_{x}^{\star} \left (\det \r p_{!}(\mathcal E\otimes
q^{\star}F) \right )^{-1} =\left ( \det \r p_{!}\left((1\times t_{rx})^{\star} \mathcal E \otimes q^{\star}\left(t_{x}^{\star}\Lambda^{-1}\otimes 
\Lambda\otimes F\right)\right) \right )^{-1}.$$
The first Chern class can then be computed by Grothendieck-Riemann-Roch. The answer does not depend on the point $x\in A$ since the maps $(1\times
t_{rx})^{\star}$ and $t_{x}^{\star}$ act as the cohomological restriction associated with $K_{v}\times A \hookrightarrow \mv \times A$, and as the identity 
on the cohomology of $A$, respectively. 
When a universal family does not exist, one can use a quasi-universal family instead. 

The above argument shows that $\left(\Phi^{+}\right)^{\star}\Theta_w$ should be of the form $\Theta_{w}\boxtimes \mathcal L^{+}$ for some line bundle 
$\mathcal
L^{+}$ coming from $A$. We can express this line bundle explicitly as follows. Write $$m=m_{1}:A\times A\to A, (a,b)\to a+b$$ for the addition map, and
consider the morphism $$m_{r}:A\times A\to A, (a, b) \to a+rb.$$ Then, $$m_{r}=m\circ (1, r).$$ Letting $p_{1}, p_{2}$ be the two projections, we have
$$\mathcal L^{+}=\left ( \det \r p_{2!}\left(m_{r}^{\star}E \otimes m^{\star}\Lambda^{-1} \otimes p_{1}^{\star}(\Lambda\otimes F)\right)\right )^{-1}.$$ 
Letting $\lambda = c_{1}(\Lambda),$ we get by
Grothendieck-Riemann-Roch, $$c_{1}(\mathcal L^{+})= - p_{2!}\left[m_{r}^{\star}v\cdot m^{\star}e^{-\lambda} \cdot
p_{1}^{\star}(e^{\lambda}w)\right]_{(3)}.$$ Expanding each of the terms, we obtain $$c_{1}(\mathcal L^{+})= - p_{2!} 
\left[\left(r+m_{r}^{\star}\lambda+\chi\,
m_{r}^{\star}\omega\right)\cdot \left(1-m^{\star}\lambda + \frac{\lambda^{2}}{2}m^{\star}\omega\right) \right.$$ $$\left. \cdot
p_{1}^{\star}\left(r'+\left(r'\lambda+\lambda'\right)+\left(\chi'+\lambda \lambda'+ r'\frac{\lambda^{2}}{2}\right)\omega\right)\right]_{(3)}.$$ The precise
evaluation of this product relies on the following intersections $$p_{2!} (m^{\star}\lambda\cdot p_{1}^{\star}\omega )=\lambda,\,\,\,\, p_{2!}
(m_{r}^{\star}\lambda\cdot p_{1}^{\star}\omega )=r^{2}\lambda,$$ $$p_{2!} (m_{r}^{\star}\omega\cdot m^{\star}\lambda)=(r-1)^{2}\lambda,\,\,\, p_{2!}
(m_{r}^{\star}\lambda \cdot m^{\star}\omega)=(r-1)^{2}\lambda,$$ $$p_{2!}(m^{\star} \omega \cdot p_{1}^{\star}\alpha)=\alpha,\,\,\,
p_{2!}(m_{r}^{\star}\omega\cdot p_{1}^{\star}\alpha)=r^{2}\alpha, \, \, \text{for any} \, \, \alpha \in H^2(A).$$ The last pair of intersections is to be 
used for the class $\alpha=r'\lambda+\lambda'.$
The formulas above are easily justified either by explicit computations in coordinates, or directly, by interpreting geometrically the intersections
involved. For instance, the third pushforward $p_{2!} (m_{r}^{\star}\omega\cdot m^{\star}\lambda)$ is computed as the image under $p_{2}$ of the cycle 
$$\left\{(a, b), \,\, a+rb=0, \,\,\, a+b\in
\lambda\right\}\hookrightarrow A\times A.$$ This pushforward can be identified with $(r-1)^{\star}\lambda=(r-1)^{2}\lambda.$

The value of the Chern class is obtained immediately from the previous intersections and a last one calculated by the 
Lemma below. Equation \eqref{relation} has to be used to bring the answer in the form claimed by Proposition \ref{split}. \vskip.1in

\begin {lemma} For any $\lambda, \alpha \in H^{2} (A)$, we have $$p_{2!}(m_{r}^{\star}\lambda\cdot m^{\star}\lambda \cdot p_{1}^{\star}\alpha)=(r-1)^{2}
\left(\int_{A}\alpha\lambda\right) \cdot \lambda + r\lambda^{2} \cdot \alpha .$$ \end {lemma}

\proof First, note the isomorphism $$m^{\star}\Lambda\cong p_{1}^{\star}\Lambda\otimes p_{2}^{\star}\Lambda \otimes (1 \times \Phi_{\Lambda})^{\star}\mathcal P.$$ This shows
that \begin{equation}\label{mstar}m^{\star}\lambda = p_{1}^{\star}\lambda+p_{2}^{\star}\lambda+(1\times \Phi_{\Lambda})^{\star}c_{1}(\mathcal P), \text { and
}\end{equation} $$m_{r}^{\star}\lambda = (1\times r)^{\star}m^{\star}\lambda=p_{1}^{\star}\lambda+r^{2}p_{2}^{\star}\lambda+r \cdot(1\times
\Phi_{\Lambda})^{\star}c_{1}(\mathcal P).$$ It follows that $$p_{2!}(m_{r}^{\star}\lambda\cdot m^{\star}\lambda \cdot
p_{1}^{\star}\alpha)=(r^{2}+1)\left(\int_{A}\alpha\lambda\right) \cdot \lambda +r \cdot p_{2!} \left(p_{1}^{\star}\alpha \cdot (1\times 
\Phi_{\Lambda})^{\star}
c_{1}(\mathcal P)^{2}\right)$$ $$= (r^{2}+1)\left(\int_{A}\alpha\lambda\right) \cdot \lambda + 2r \cdot \Phi_{\Lambda}^{\star} 
\left\{p_{2!}\left(p_{1}^{\star}\alpha \cdot
\frac{c_{1}(\mathcal P)^{2}}{2}\right)\right\}.$$ We will prove \begin{equation}\label{fmp}\Phi_{\Lambda}^{\star} \left\{p_{2!}\left(p_{1}^{\star}\alpha
\cdot \frac{c_{1}(\mathcal P)^{2}}{2}\right)\right\}=-\left(\int_{A}\alpha \lambda \right)\cdot \lambda +\frac{\lambda^{2}}{2}\cdot \alpha.\end{equation}

This follows by a computation in coordinates. Explicitly, let us write $A=V/\Gamma$. We regard $V$ as a four-dimensional real vector space. The dual
abelian variety has as underlying {\it real} manifold the torus $V^{\vee}/\Gamma^{\vee}$, where $V^{\vee}$ stands for the {\it real} dual of $V$. Pick a
basis $f_{1}, f_{2}, f_{3}, f_{4}$ for $V$, which is symplectic for $\Lambda$. This means that in the dual basis, $$\lambda=c_{1}(\Lambda)=d \cdot f_{1}^{\vee}\wedge
f_{2}^{\vee}+ e \cdot f_{3}^{\vee}\wedge f_{4}^{\vee}\in \Lambda^{2}V^{\vee},$$ for some (integers) $d$ and $e$. Moreover, the Chern class of the Poincar\'{e} line bundle on
$A\times \widehat A$ takes the form \begin{equation}\label{poincare}c_{1}(\mathcal P)=f_{1}^{\vee}\wedge f_{1}+f_{2}^{\vee}\wedge f_{2}+f_{3}^{\vee}\wedge
f_{3}+f_{4}^{\vee}\wedge f_{4}.\end{equation} To prove \eqref{fmp}, it suffices to assume that $$\alpha =f_{1}^{\vee} \wedge f_{2}^{\vee}, \text { or 
}
\alpha=f_{1}^{\vee} \wedge f_{3}^{\vee}.$$ Let us consider only the first case, the second being similar. Then, $$p_{2!}\left(p_{1}^{\star}\alpha \cdot 
\frac{c_{1}(\mathcal
P)^{2}}{2}\right)=-f_{3} \wedge f_{4}.$$ The discussion in \cite {LB}, chapter $2$, and in particular Lemma $4.5$ therein, shows that the map
$$\Phi_{\Lambda}^{\star}:H^{1}(\widehat A, \mathbb R)\cong V \to H^{1}(A, \mathbb R)\cong V^{\vee}$$ is induced by the contraction of the first Chern class
$c_{1}(\Lambda)$.  It follows that $$\Phi_{\Lambda}^{\star}\left\{ p_{2!}\left(p_{1}^{\star}\alpha \cdot \frac{c_{1}(\mathcal
P)^{2}}{2}\right)\right\}=-\Phi_{\Lambda}^{\star}f_{3} \wedge \Phi_{\Lambda}^{\star} f_{4}= - e^{2}\, f_{3}^{\vee}\wedge f_{4}^{\vee}.$$ But this is also
the result on the right hand side of \eqref{fmp}:$$-\left(\int_{A}\alpha \lambda \right)\cdot \lambda +\frac{\lambda^{2}}{2}\cdot \alpha=-e\cdot\lambda + de
\cdot \alpha = -e^{2}\, f_{3}^{\vee}\wedge f_{4}^{\vee}.$$
 \vskip.1in {\it Proof of Theorem \ref {main}.} When $d_{v}\neq 0$, Theorem \ref {main} follows immediately from Propositions \ref{kv} and \ref {split}, and 
Lemma \ref{degree1}. Indeed, we
have $$\chi(\mv, \Theta_{w})=\frac{1}{d_{v}^{4}}\chi(\left(\Phi^{+}\right)^{\star}\Theta_{w})=\frac{1}{d_{v}^{4}}\chi(K_{v}, \Theta_{w})\,\chi (A, \mathcal
L^{+})$$ $$=\frac{1}{d_{v}^{4}} \cdot\frac{d_{v}^{2}}{d_{v}+d_{w}}\binom{d_{v}+d_{w}}{d_{v}}\cdot \frac{\left(d_{v} c_{1}(v\otimes
w)\right)^{2}}{2}=\frac{1}{2}\frac{c_{1}(v\otimes w)^{2}}{d_{v}+d_{w}} \binom {d_{v}+d_{w}}{d_{v}}.$$

When $d_{v}=0$, the Theorem is equivalent to the equality $$\chi(\mv, \Theta_{w})=r^{2}.$$ It suffices to explain that the moduli space $\mv$ consists of
$r^{2}$ smooth points. By work of Mukai, it is known that $\uv$ is an abelian surface. In fact, fixing $E\in \mv$, we have an isogeny $$A\to \uv,
\,\,\,x\to t_{x}^{\star}E,$$ whose kernel is the group $$K(E)=\left\{x \in A\text { such that } t_{x}^{\star}E\cong E\right\}.$$ Note that we may need to replace $E$ by a twist $E\otimes H^{\otimes n}$ to ensure that $K(E)$ is finite. In this case, $K(E)$ has $\chi^{2}$ elements. This is a result of Mukai \cite {mukai2}; to apply it, we need to observe that $E$ is a simple semi-homogeneous sheaf. Restricting to sheaves with determinant $\Lambda$, we see that $$\mv\cong K(\Lambda)/K(E),$$ has length $\frac{\chi(\Lambda)^{2}}{\chi^{2}}=r^{2}.$

\section {Sheaves with fixed determinant of the Fourier-Mukai transform.} This section is devoted to the proof of Theorem \ref {two}. It is possible to 
deduce this result
from Theorem \ref {main} when the Picard number of $A$ is $1$, by explicitly studying how the relevant moduli spaces and Theta divisors change under the
Fourier-Mukai transform \cite {Y}\cite {mac}. However, the following proof is simpler, covers all cases, and it is in the spirit of this paper. Note
that the cohomological computation below may be regarded as the Fourier-Mukai 'dual' of last section's calculations.

We will crucially make use of the diagram \begin{center} $\xymatrix{K_{v}\times \widehat {A} \ar[r]^{\Phi^{-}} \ar[d]^{p} & \mvv
\ar[d]^{\alpha^{+}} \\ {\widehat A}\ar[r]^{\Psi^{-}} & {\widehat A}}\label{diag1}.$ \end{center} The upper horizontal morphism $\Phi^{-}$ is defined as
$$\Phi^{-}(E, y)=t_{\Phi_{\widehat \Lambda}(y)}^{\star}E\otimes y^{\chi}.$$ To check that $\Phi^{-}$ is well defined, we compute 
\begin{eqnarray*}\nonumber\det \fm
\left(\Phi^{-}(E, y)\right)&=&\det \fm \left (t_{\Phi_{\widehat \Lambda}(y)}^{\star}E\otimes y^{\chi}\right)=\det \left(t_{\chi y}^{\star}\, \fm (E)\otimes
\Phi_{\widehat\Lambda}(y)^{-1}\right) \\&=&t_{\chi y}^{\star} \det \fm (E) \otimes \Phi_{\widehat \Lambda}(y)^{-\chi}=t_{\chi y}^{\star} \widehat \Lambda \otimes
\Phi_{\widehat \Lambda}(y)^{-\chi}\\ &=&\widehat\Lambda\otimes \Phi_{\widehat \Lambda}(\chi y)\otimes \Phi_{\widehat \Lambda}(y)^{-\chi}=\widehat \Lambda.\end
{eqnarray*} The next two results are the versions of Lemma \ref{multlemma} and Proposition \ref{split} suitable to the present context.

\begin {lemma} The morphism $\Psi^{-}$ is the multiplication by $-d_{v}$ on the abelian variety $\widehat A$.
\label{mult2}
\end {lemma}

\proof Using \eqref{fmpp} and \eqref{phis}, we compute $$\alpha^{+}\circ \Phi^{-}(E, y)=\det \left(t_{\Phi_{\widehat \Lambda}(y)}^{\star}E\otimes
y^{\chi}\right)\otimes \Lambda^{-1}=t_{\Phi_{\widehat \Lambda}(y)}^{\star}\Lambda \otimes y^{r\chi}\otimes \Lambda^{-1}$$
$$=\Phi_{\Lambda}\left(\Phi_{\widehat \Lambda}(y) \right)\otimes y^{r\chi}=(- \chi(\Lambda)+r\chi)y=-d_{v}y.$$

\begin {proposition}\label{split1} We have $$\left(\Phi^{-}\right)^{\star} \Theta_{w}\cong \Theta_{w}\boxtimes \mathcal L^{-},$$ where $$c_{1}(\mathcal
L^{-})= - d_{v}c_{1}(\hat v \otimes \hat w).$$ \end {proposition} \proof The proof of this result parallels that of Proposition \ref{split}. It suffices to 
show that 
the
line bundle $\mathcal L^{-}$ corresponding to the divisor $$\left\{y\in \widehat A,\,\,\text {with } H^{0} \left(t^{\star}_{\widehat
\Lambda(y)}(E)\otimes y^{\chi}\otimes F\right)\neq 0\right\}$$ has the first Chern class given by the Proposition. Note that $$\mathcal L^{-}= \left (\det
p_{2!}\left(f^{\star}E\otimes p_{1}^{\star}F\otimes \mathcal P^{\chi}\right) \right)^{-1},$$ where $$f:A\times \widehat A\to A\times A{\to}A$$ denotes the 
composition
\begin{equation}\label{comp}f=m\circ (1 \times \Phi_{\widehat \Lambda}),\,\, (x, y)\to x+\Phi_{\widehat \Lambda}(y).\end{equation} By Riemann-Roch, we
compute $$c_{1}(\mathcal L^{-})= - p_{2!}\left[\left(r+f^{\star}\lambda+\chi f^{\star}\omega\right)\cdot
\left(r'+p_{1}^{\star}\lambda'+\chi'p_{1}^{\star}\omega\right)\cdot \left(1+\chi c_{1}(\mathcal P)+\chi^{2}\frac{c_{1}(\mathcal P)^{2}}{2}\right)
\right]_{(3)}.$$ The following observations allow for the explicit evaluation of the expression above: $$p_{2!}\left(\frac{c_{1}(\mathcal P)^{2}}{2}\cdot
p_{1}^{\star}\lambda'\right)=\widehat \lambda',\,\, p_{2!}\left(\frac{c_{1}(\mathcal P)^{2}}{2}\cdot f^{\star}\lambda\right)=\widehat \lambda,$$
\begin{equation}\label{a}p_{2!}(f^{\star}\omega\cdot p_{1}^{\star}\lambda')=\frac{\lambda^{2}}{2}\cdot \widehat\lambda'-(\lambda\cdot \lambda')\cdot \widehat
\lambda,\end {equation} \begin {equation}\label{b} p_{2!}\left(f^{\star}\lambda \cdot p_{1}^{\star}\omega\right)=- \frac{\lambda^{2}}{2}\cdot\widehat\lambda,
\end {equation} \begin {equation}\label{c} p_{2!}\left(f^{\star}\omega \cdot c_{1}(\mathcal P)\right)=-2\widehat \lambda, \end {equation} \begin
{equation}\label{d} p_{2!}\left(f^{\star}\lambda \cdot p_{1}^{\star}\lambda'\cdot c_{1}(\mathcal P)\right)=-\lambda^{2}\cdot \widehat \lambda'. \end
{equation} The Proposition follows by substitution, also making straightforward use of the orthogonality constraint $$r\chi'+\lambda\cdot \lambda'+r'\chi=0.$$

It remains to explain the four numbered equations claimed above. Let us first consider \eqref{a}. Interpreting the pushforward geometrically, and recalling 
the
definiton of $f$ in \eqref{comp}, we find that $$p_{2!}(f^{\star}\omega\cdot p_{1}^{\star}\lambda')=(-\Phi_{\widehat \Lambda})^{\star}\lambda'=\Phi_{\widehat
\Lambda}^{\star}\lambda'=\frac{\lambda^{2}}{2}\cdot \widehat\lambda'-\left(\int_{A}\lambda\cdot \lambda'\right)\cdot \widehat \lambda.$$ The dual of the last 
equality
was verified in \eqref{fmp}. The case at hand is a corollary of what we have already shown there, using the fact that the
Fourier-Mukai transform is an isometry. Equation \eqref{b} is very similar. To prove it, we observe that $f$ restricts to $\Phi_{\widehat\Lambda}$ on
$\{0\}\times\widehat A$, hence $$p_{2!}\left(f^{\star}\lambda \cdot p_{1}^{\star}\omega\right)=\Phi_{\widehat
\Lambda}^{\star}\lambda=-\frac{\lambda^{2}}{2}\cdot \widehat {\lambda}.$$

In turn, \eqref{c} follows by a computation in local coordinates. First, pick a basis $f_{1}, f_{2}, f_{3}, f_{4}$ for $V$ such that $$c_{1}(\mathcal 
P)=f_{1}^{\vee}\wedge
f_{1}+f_{2}^{\vee}\wedge f_{2}+f_{3}^{\vee}\wedge f_{3}+f_{4}^{\vee}\wedge f_{4}.$$ From the definition of $f$ in \eqref{comp}, we calculate
$$p_{2!}(f^{\star}\omega \cdot c_{1}(\mathcal P))= p_{2!}\left((1\times\Phi_{\widehat \Lambda})^{\star}m^{\star}\omega \cdot c_{1}(\mathcal P)\right)= $$
$$=-p_{2!}\left((1\times \Phi_{\widehat \Lambda})^{\star} \left(\sum_{j=1}^{4}\text{PD}(f_{j}^{\vee}) \wedge f_{j}^{\vee}\right)\cdot \left(\sum_{j=1}^{4} 
f_{j}^{\vee}\wedge
f_{j}\right)\right)=\sum_{j=1}^{4}\Phi_{\widehat \Lambda}^{\star}f_{j}^{\vee} \wedge f_{j}.$$ Taking 
$$\lambda = d
\cdot f_{1}^{\vee}\wedge f_{2}^{\vee}+e\cdot f_{3}^{\vee}\wedge f_{4}^{\vee},$$ this last expression is 
\begin{equation}\label{fmtl} 
2d \cdot f_{3}\wedge f_{4}+2e\cdot f_{1}\wedge f_{2} = -2 \widehat{\lambda},\end{equation}
confirming \eqref{c}.

Finally, for \eqref{d}, we observe that \begin{eqnarray*}p_{2!}\left(f^{\star}\lambda \cdot p_{1}^{\star}\lambda'\cdot c_{1}(\mathcal
P)\right)&=&p_{2!}\left((1\times \Phi_{\widehat\Lambda})^{\star}m^{\star}\omega\cdot p_{1}^{\star}\lambda'\cdot c_{1}(\mathcal P)\right) \\
&=&p_{2!}\left((1\times \Phi_{\widehat\Lambda})^{\star}(1\times \Phi_{\Lambda})^{\star}c_{1}(\mathcal P)\cdot p_{1}^{\star}\lambda'\cdot c_{1}(\mathcal P)\right)
\\ &=& p_{2!}\left(\left(1\times \left(- \chi(\Lambda\right)\right)\right)^{\star}c_{1}(\mathcal P)\cdot p_{1}^{\star}\lambda'\cdot c_{1}(\mathcal P))\\ &=&-\chi(\Lambda)\cdot
p_{2!}(c_{1}(\mathcal P)^{2}\cdot p_{1}^{\star}\lambda')=\lambda^{2}\cdot \widehat \lambda'.\end {eqnarray*} The first line follows by the definition of $f$
in \eqref{comp}, the second uses \eqref{mstar}, while the third uses \eqref{phis}.

\vskip.1in {\it Proof of Theorem \ref{two}.} As before, when $d_{v}\neq 0$, the Theorem follows immediately from Propositions \ref{kv} and \ref {split1}, and 
Lemma \ref{mult2}. Using
the cartesian diagram, we compute $$\chi(\mvv, \Theta_{w})=\frac{1}{d_{v}^{4}}\chi(\left(\Phi^{-}\right)^{\star}\Theta_{w})=\frac{1}{d_{v}^{4}}\chi(K_{v},
\Theta_{w})\,\chi (A, \mathcal L^{-})$$ $$=\frac{1}{d_{v}^{4}} \cdot\frac{d_{v}^{2}}{d_{v}+d_{w}}\binom{d_{v}+d_{w}}{d_{v}}\cdot \frac{\left(d_{v} c_{1}(\hat
v\otimes \hat w)\right)^{2}}{2}=\frac{1}{2}\frac{c_{1}(\hat v\otimes \hat w)^{2}}{d_{v}+d_{w}} \binom {d_{v}+d_{w}}{d_{v}}.$$

When $d_{v}=0$, we observe that $\mvv$ consists of $\chi^{2}$ smooth points. First, for any sheaf $E$ in the moduli space $\uv$, consider  
the isogeny $$\widehat A\to \uv,\, y\mapsto E\otimes \mathcal P_{y}.$$The kernel $$\Sigma(E)=\left\{y \in \widehat A \text { such that } E\otimes \mathcal P_{y} \cong E\right\}$$ has
length $r^{2}$, cf. \cite {mukai2} (twisting by powers of $H$ may be necessary). Note that the points in $\mvv$ have the property $$\det \fm \left(E\otimes \mathcal P_{y}\right)\otimes (\det \fm E)^{-1}\cong
t_{y}^{\star}\widehat \Lambda \otimes \widehat \Lambda^{-1}\cong \mathcal O.$$ Therefore, $$\mvv=K(\widehat \Lambda)/\Sigma(E)$$ has length $\chi(\widehat
\Lambda)^{2}/r^{2}=\chi^{2}$, as claimed.

\section {Sheaves with arbitrary determinant.} This last section contains the proof of Theorem \ref{three}. The Euler characteristic on $K_{v}$ was calculated
in Proposition \ref{kv}. To compute the one on $\uw$, we use the diagram \begin{center} $\xymatrix{K_{w}\times A \times \widehat
A\ar[r]^{\Phi} \ar[d]^{p} & \uw \ar[d]^{\mathsf a_{w}} \\ A\times \widehat A \ar[r]^{\Psi} & A\times \widehat A}.$ \end{center} Here,  
$\Phi:K_{w}\times
A\times \widehat A\to \uw$ is defined as $$\Phi(E, x, y)=t_{x}^{\star}E \otimes y.$$ Using \eqref{fmpp1} and \eqref{fmpp}, Yoshioka proved in detail that $$\Psi(x, y)=(-\chi' x+ \Phi_{\widehat \Lambda'}(y), \Phi_{\Lambda'}(x)+r'y),$$ which has degree
$d_{w}^{4}$ \cite {Y}.

\begin {proposition}\label{split2} We have $$\Phi^{\star}\Theta_{v}=\Theta_{v}\boxtimes \mathcal L$$ where $\mathcal L$ is a line bundle on $A\times \widehat 
A$
with $$\chi(\mathcal L)=d_{v}^{2}d_{w}^{2}.$$ \end {proposition}

\proof It suffices to compute the Euler characteristic of the line bundle $\mathcal L$ corresponding to the divisor $$\{(x, y)\in A\times \widehat A, \text {
such that } H^{0}(t_{x}^{\star}E\otimes y\otimes F)\neq 0\}.$$ In other words $$\mathcal L=\left ( \det p_{23!}\left(m_{12}^{\star}E\otimes 
p_{13}^{\star}\mathcal
P\otimes p_{1}^{\star}F\right) \right)^{-1},$$ where the $p$'s denote the projections on the corresponding factors of $A\times A\times \widehat A,$ while 
$$m_{12}:A\times
A\times \widehat A\to A$$ is the addition on the first two factors.  Keeping the previous notations, $$c_{1}(\mathcal
L)= - p_{23!}\left[\left(r+m_{12}^{\star}\lambda + \chi m_{12}^{\star}\omega \right)\cdot \left(1+p_{13}^{\star}c_{1}(\mathcal
P)+\frac{p_{13}^{\star}c_{1}(\mathcal P)^{2}}{2}\right) \cdot \left(r'+p_{1}^{\star}\lambda'+\chi' p_{1}^{\star}\omega\right) \right]_{(3)}.$$ Expanding, we
easily obtain $$ - c_{1}(\mathcal L)=(\chi\lambda'+\chi'\lambda)+(r\widehat \lambda'+r'\widehat \lambda)-r\chi' c_{1}(\mathcal
P)+p_{23!}\left(m_{12}^{\star}\lambda\cdot p_{13}^{\star}c_{1}(\mathcal P)\cdot p_{1}^{\star}\lambda'\right).$$

We claim that $$\chi(\mathcal L)=\frac{c_{1}(\mathcal L)^{4}}{4!}=d_{v}^{2}d_{w}^{2}.$$ The computation makes use of the fact that the Picard 
number of $A$ is
$1$, so we may assume that either $\lambda'=0$, or otherwise that $$\lambda=a \lambda'$$ for some constant $a$. In the first case, we have
$$c_{1}(\mathcal L)=- \chi'\lambda - r'\widehat \lambda + r\chi' c_{1}(\mathcal P).$$ To prove the claim, we first note that \begin{equation}\label{llp} \lambda\cdot 
\widehat\lambda\cdot
\frac{c_{1}(\mathcal P)^{2}}{2}=\lambda^{2}\end{equation} This follows easily by a computation in local coordinates. Indeed, writing
\begin{equation}\label{local}\lambda=d \,f_{1}^{\vee}\wedge f_{2}^{\vee}+ e \, f_{3}^{\vee}\wedge f_{4}^{\vee},\end{equation} and recalling that $\widehat 
\lambda$ and $c_{1}(\mathcal
P)$ have the form \eqref{fmtl} and \eqref{poincare} respectively, we calculate $$\lambda\cdot \widehat \lambda \cdot \frac{c_{1}(\mathcal
P)^{2}}{2}=2de=\lambda^{2}.$$ With \eqref{llp} understood, and using the fact that the Fourier-Mukai is an isometry, we obtain $$\frac{c_{1}(\mathcal
L)^{4}}{4!}=\frac{1}{4!}(\chi' \lambda+r' \widehat \lambda-r\chi'c_{1}(\mathcal
P))^{4}=\frac{r'^{2}\chi'^{2}(\lambda^{2})^{2}}{4}+r'\chi'(r\chi')^{2}\lambda^{2}+(r\chi')^{4}$$
$$=\left(\frac{r'\chi'\lambda^{2}}{2}+(r\chi')^{2}\right)^{2}=\left[r'\chi'\left(\frac{\lambda^{2}}{2}-r\chi\right)\right]^{2}=d_{v}^{2}d_{w}^{2}.$$ The
penultimate equality made use of the fact that $r\chi'+r'\chi=0.$

Finally, the more general second case $\lambda = a \lambda'$ is similar. Using \eqref{mstar}, we get $$p_{23!}\left(m_{12}^{\star}\lambda\cdot
p_{13}^{\star}c_{1}(\mathcal P)\cdot p_{1}^{\star}\lambda'\right)=(\Phi_{\Lambda}\times 1)^{\star}q_{23!}\left(q_{12}^{\star}c_{1}(\mathcal P)\cdot
q_{13}^{\star}c_{1}(\mathcal P)\cdot q_{1}^{\star}\lambda'\right),$$ with the $q$'s standing for the projections of the factors of $A\times \widehat A\times
\widehat A.$ In turn, we claim that \begin{equation}\label{bl}(\Phi_{\Lambda}\times 1)^{\star}q_{23!}\left(q_{12}^{\star}c_{1}(\mathcal P)\cdot
q_{13}^{\star}c_{1}(\mathcal P)\cdot q_{1}^{\star}\lambda\right)=-\frac{\lambda^{2}}{2}c_{1}(\mathcal P).\end{equation}Again, this is easiest to check in
local coordinates. Assuming that \eqref{local} holds, we have $$q_{23!}\left(q_{12}^{\star}c_{1}(\mathcal P)\cdot q_{13}^{\star}c_{1}(\mathcal P)\cdot
q_{1}^{\star}\lambda\right)=-d f_{3} \otimes f_{4} +d f_{4}\otimes f_{3} - ef_{1}\otimes f_{2}+ef_{2}\otimes
f_{1}.$$ Hence, after pullback by $\Phi_{\Lambda}\times 1$, the left hand side of \eqref{bl} becomes $$-de \,c_{1}(\mathcal
P)=-\frac{\lambda^{2}}{2}c_{1}(\mathcal P),$$ as claimed. Putting things together, we obtain $$c_{1}(\mathcal L)= - (\chi'a+\chi)\lambda' - (r'a +r)\widehat
\lambda' +\left(r\chi'+\frac{a \lambda'^{2}}{2}\right)c_{1}(\mathcal P).$$ The same type of calculation as the one done above yields the answer
$$\chi(\mathcal L)=\frac{c_{1}(\mathcal
L)^{4}}{4!}=\left[\frac{(\chi'a+\chi)(r'a+r)\lambda'^{2}}{2}+\left(r\chi'+\frac{a{\lambda'}^{2}}{2}\right)^{2}\right]^{2}.$$ To conclude the proof, it
remains to observe that the expression in square brackets can be equated with $$-\left(\frac{a^{2}\lambda'^{2}}{2}-r\chi
\right)\left(\frac{{\lambda'}^{2}}{2}-r'\chi'\right)=-d_{v}d_{w},$$ so that $$\chi(\mathcal L)=(d_{v}d_{w})^{2}.$$ The latter algebraic manipulation will be
left to the reader, who may wish to use the fact that $$a \lambda'^{2}+r\chi'+r'\chi=0.$$

It is very likely that the Lemma holds true for arbitrary abelian surfaces, without any restrictions on the N\'{e}ron-Severi group, but the computation seems 
to
be more involved.

\vskip.1in {\it Proof of Theorem \ref{three}.} We compute $$\chi(\uw, \Theta_{v})=\frac{1}{d_{w}^{4}}\chi(K_{w}, \Theta_{v})\chi (A \times \widehat A,
\mathcal L)=\frac{1}{d_{w}^{4}} \cdot \frac{d_{w}^{2}}{d_{v}+d_{w}}\binom{d_{v}+d_{w}}{d_{v}}\cdot \left(d_{v}^{2}d_{w}^{2}\right)$$
$$=\frac{d_{v}^{2}}{d_{v}+d_{w}}\binom{d_{v}+d_{w}}{d_{v}}=\chi(K_{v}, \tw).$$ The case $d_{w}=0$ requires, as usual, special care. We need to show
$$\chi(\uw, \tv)=d_{v}.$$ Using the degree $\chi'^{2}$ isogeny: $$\pi:A\to \uw, x\to t_{x}^{\star}F,$$ where $F$ is a semi-homogeneous sheaf of Mukai
vector $w$, we have $$\pi^{\star}\Theta_{v}=\left ( \det p_{!}(m^{\star}F\otimes q^{\star}E) \right )^{-1},$$ with $p$, $q$ and $m$ standing for the
projection and addition morphism. Then $$c_{1}(\pi^{\star}\Theta_{v})= - \chi' \lambda - \chi \lambda'.$$ We obtain $$\chi(\uw,
\tv)=\frac{1}{\chi'^{2}}\chi(A, \pi^{\star}\tv)=\frac{1}{2\chi'^{2}}(\chi'\lambda+\chi\lambda')^{2}=\frac{1}{2\chi'^{2}}\left(\chi'^{2}\lambda^{2}+\chi^{2}
\lambda'^{2}-2\chi\chi'(r\chi'+r'\chi)\right)$$
$$=\frac{1}{\chi'^{2}}(\chi^{2}d_{w}+\chi'^{2}d_{v})=d_{v}.$$ This completes the proof of the Theorem.


\begin{thebibliography}{1}

\bibitem [B]{B}

P. Belkale, {\it The strange duality conjecture for generic curves}, J. Amer Math Soc,  to appear.

\bibitem [D]{D}

G. D\v{a}nil\v{a}, {\it R\'{e}sultats sur la conjecture de dualit\'{e} \'{e}trange sur le plan projectif}, Bull. Soc. Math. France, 130 (2002), no. 
1, 1-33.

\bibitem [DN]{dn}

J. M. Dr\'{e}zet, M.S. Narasimhan, {\it Groupe de Picard des vari\'{e}t\'{e}s de modules de fibr\'{e}s semi-stables sur les courbes alg\'{e}briques}, Invent. Math. 97 (1989), no. 1, 53-94.

\bibitem [EGL]{EGL}

G. Ellingsrud, L. G\"{o}ttsche, M. Lehn, {\it On the cobordism class of the Hilbert scheme of a surface,} J. Algebraic Geom. 10 (2001), 81-100.

\bibitem [GNY]{GNY}

L. G\"{o}ttsche, H. Nakajima, K. Yoshioka, {\it K-theoretic Donaldson invariants via instanton counting}, math-AG/061194.

\bibitem [HT]{HT}

B. Hassett, Y. Tschinkel, {\it Rational curves on holomorphic symplectic fourfolds},   Geom. Funct. Anal. 11 (2001), 1201-1228. 

\bibitem [H]{H1}
 
D.  Huybrechts, {\it Compact Hyperk\"{a}hler manifolds}, Calabi-Yau manifolds and related geometries, Lectures from the Summer School held in Nordfjordeid, 
June 2001 Springer-Verlag, Berlin, 2003, 161-225.


\bibitem [LB]{LB}

H. Lange, C. Birkenhake, {\it Complex abelian varieties}, Springer-Verlag, Berlin-New York, 1992. 

\bibitem [LP]{LP}

J. Le Potier, {\it Fibr\'{e} d\'{e}terminant et courbes de saut sur les surfaces alg\'{e}briques,} Complex Projective Geometry, 
213-Ð240, London Math. Soc. Lecture Note Ser., 179, Cambridge Univ. Press, Cambridge, 1992.

\bibitem [Li]{Li}

J. Li, {\it Algebraic geometric interpretation of Donaldson's polynomial invariants,} J. Differential Geom. 37 (1993), 417-466.

\bibitem [MO1]{MO}

A. Marian, D. Oprea, {\it The level-rank duality for non-abelian theta functions}, Invent. Math. 168 (2007), 225-247.

\bibitem[MO2] {mo2}

A. Marian, D. Oprea, {\it A brief tour of theta dualities on moduli spaces of sheaves}, preprint (2007).

\bibitem [OG]{OG}

K. O'Grady, {\it Involutions and linear systems on holomorphic symplectic manifolds,} Geom. Funct. Anal. 15 (2005), no. 6, 1223-1274.

\bibitem [Ma]{mac}

A. Maciocia, {\it The determinant line bundle over moduli spaces of instantons on abelian surfaces}, Math. Z. 217 (1994), no. 2, 317-333.

\bibitem [Muk1]{Muk}

S. Mukai, {\it Duality between $D(X)$ and $D(\hat X)$ with its application to Picard sheaves},  Nagoya Math. J. 81 (1981), 153-175.

\bibitem [Muk2]{mukai2}

S. Mukai, {\it Semi-homogeneous vector bundles on an Abelian variety}, J. Math. Kyoto
Univ. 18 (1978), 239-272.

\bibitem [Muk3] {mukai3}

S. Mukai, {\it On the moduli space of bundles on $K3$ surfaces I}, Vector bundles on algebraic varieties, Tata Institute for Fundamental Research Studies in
Mathematics, no. 11 (1987), 341-413.

\bibitem [M]{M}

D. Mumford, {\it Abelian varieties}, Tata Institute of Fundamental Research Studies in Mathematics, Oxford University Press, London, 1970.

\bibitem [Y1]{Y}

K. Yoshioka, {\it Moduli spaces of stable sheaves on abelian surfaces}, Math. Ann. 321 (2001), 817-884.

\bibitem [Y2] {Y3}

K. Yoshioka, {\it Some notes on the moduli of stable sheaves on elliptic surfaces}, Nagoya Math. J. 154 (1999), 73-102.


\end{thebibliography}
\end{document}